\documentclass[11pt,a4paper]{amsart}
\usepackage{amsmath,amssymb,
cleveref,enumitem,graphicx
}
\renewcommand{\uppercasenonmath}[1]{}

\newcommand{\rB}{\rotatebox[origin=c]}
\newcommand{\lA}{\mathrel{\rB{180}{\circlearrowright}}}
\newcommand{\rA}{\mathrel{\rB{180}{\circlearrowleft}}}
\setlist[enumerate]{label*=(\arabic*)}
\theoremstyle{definition}
\newtheorem{dfn}{Definition}[section]
\newtheorem{rmk}[dfn]{Remark}
\newtheorem{exm}[dfn]{Example}
\newtheorem{fct}[dfn]{Fact}

\newtheorem{lmm}[dfn]{Lemma}
\newtheorem{crr}[dfn]{Corollary}
\newtheorem{prp}[dfn]{Proposition}
\newtheorem{thm}[dfn]{Theorem}
\newtheorem*{mth}{Main Theorem}

\begin{document}

\title[Local Langlands correspondence for covering groups of tori]%
{Local Langlands correspondence for covering groups of tori, 
and the packet-indexing groups}
\author[Y.\ Nakata]{Yuki Nakata}
\email{nakata.yuuki.42x@st.kyoto-u.ac.jp}

\address{Department of Mathematics, Graduate School of Science, 
Kyoto University}
\maketitle
\begin{abstract}
  We determine the finite group $\mathcal S$ parametrizing 
  a packet in the local Langlands correspondence 
  for a Brylinski-Deligne covering group of an algebraic torus,
  under some assumption on ramification.
  Especially, this work generalizes 
  Weissman's result on covering groups of
  tori that split over an unramified extension of the base field.
\end{abstract}
\tableofcontents

\section{Introduction}
This paper aims to determine the finite group
$\mathcal S=\mathcal S_{\widetilde T}$ whose dual parametrizes 
a packet in the local Langlands correspondence for a 
covering group $\widetilde T\to T$ of an algebraic torus.
As well as reductive algebraic groups, their covering groups
have been studied with arithmetic interest.
For instance, Weil studied the Segal-Shale-Weil representations
of the non-trivial double cover 
$\operatorname{Mp}_{2r}\to\operatorname{Sp}_{2r}$, in order to
prove an identity now called the Siegel-Weil formula~\cite{Weil_Sieg}.
Another notable example is an Eisenstein series
on a double cover of $\operatorname{GL}_r$, which gives
a Rankin-Selberg integral for the symmetric square L-function
of a cuspidal representation of $\operatorname{GL}_r$~\cite{BpG_SymS}.
Such examples motivate us to pursue the Langlands program
for covering groups.

In the familiar case of a reductive algebraic group $G$ 
over a non-archimedean local field $F$,
we know the notions of the complex dual $\widehat G$ of the group $G$
and a Langlands parameter
$\phi\colon\mathrm W_F\times\operatorname{SL}_2(\mathbb C)\to
\vphantom{G}^{\mathrm L}G=\widehat G\rtimes\mathrm W_F$.
It has been conjectured that the centralizer 
$S_\phi=\operatorname{Cent}_{\widehat G}(\operatorname{Im}\phi)$
indexes the representations in the L-packet $\Pi_\phi$,
through its finite quotient
$\mathcal S_\phi=\left.S_\phi\middle/
S_\phi^0\operatorname Z(\widehat G)^{\mathrm W_F}\right.$.
Here, $S_\phi^0$ is the connected component of the identity
in the complex reductive group $S_\phi$, and
$\operatorname Z(\widehat G)^{\mathrm W_F}$
is the subgroup consisting of fixed points in the center 
$\operatorname Z(\widehat G)\subset\widehat G$
under the action of the Weil group $\mathrm W_F$ 
of the base field $F$.
In this paper, we study an analogue of the group $\mathcal S_\phi$
for covering groups.

Though a local Langlands correspondence for a general covering group
has not been satisfactorily formulated,
one for a covering $\widetilde T\to T$ of a torus 
is proposed by Weissman~\cite{Weiss_loc_tori}.
According to his work, a certain finite abelian group
$\mathcal S=\mathcal S_{\widetilde T}$
parametrizes the representations in a packet (see \Cref{eLlc}),
as $\mathcal S_\phi$ does in the non-cover case.
He also explicitly determined the group 
$\mathcal S_{\widetilde T}$ for a cover of an unramified torus $T$, 
i.e.
a torus splitting over an unramified extension of the base field.

To describe the group $\mathcal S_{\widetilde T}$,
Weissman introduced certain subgroups $Y^\#$ and $Y^{\Gamma\#}$
of the cocharacter lattice $Y$ of $T$.
In his formulation of the local Langlands correspondence,
the cover 
$\mu_n\to\widetilde T\to T$ is assumed to be
a certain functorial one
defined by Brylinski and Deligne~\cite{Bry_Del},
which induces the groups $Y^\#$ and $Y^{\Gamma\#}$ as follows
(see \Cref{eCov}).
First, such a cover $\widetilde T$ is known to be associated with
a bilinear form $B_{\widetilde T}\colon Y\times Y\to\mathbb Z$.
Then for a subgroup $Y'\subset Y$, 
the subgroup $Y'^\#$ of $Y$ is defined by
\[Y'^\#=\{y\in Y\mid B_{\widetilde T}(y,Y')\subset n\mathbb Z\}.\]
In the present case, this subgroup $Y'$ is $Y$ itself or
the subgroup $Y^\Gamma$ of fixed points under the action
of the absolute Galois group 
$\Gamma=\Gamma_F=\operatorname{Gal}(\overline F/F)$ of $F$.

For a subgroup $Y'$ of $Y$,
let $\iota$ denote the inclusion map $Y'\hookrightarrow Y$
and the induced ones, e.g.
$\iota\colon(Y'\otimes\overline{\mathbf f}^\times)^\Gamma\to
(Y\otimes\overline{\mathbf f}^\times)^\Gamma$,
where $\overline{\mathbf f}^\times$ is the multiplicative group
of an algebraic closure of the residue field $\mathbf f$
of the valuation ring in $F$.
This notation is useful to state the following theorem by Weissman.

\begin{thm}[Weissman~\cite{Weiss_loc_tori}]
  Let $F$ be a non-archimedean local field of characteristic zero,
  $\mathbf f$ the residue field of the valuation ring in $F$, and
  $\mu_n\subset \mathbf f^\times$ the cyclic subgroup of order $n$.
  Suppose that an algebraic torus $T$ defined over $F$ 
  splits over an unramified extension $L/F$, and that
  $\mu_n\to\widetilde{T}\to T$ is a Brylinski-Deligne covering group.

  Then the group $\mathcal S=\mathcal S_{\widetilde T}$ 
  parametrizing a packet of representations of $\widetilde{T}$ is 
  written as the quotient
  \[
    \mathcal S_{\widetilde T}
    =\left.\iota((Y^{\Gamma\#}\otimes_{\mathbb Z}
    \overline{\mathbf f}^\times)^\Gamma)\middle/
    \iota((Y^{\#}\otimes_{\mathbb Z}
    \overline{\mathbf f}^\times)^\Gamma)\right..
  \]
  of subgroups of 
  $(Y\otimes_{\mathbb Z}\overline{\mathbf f}^\times)^\Gamma$.
\end{thm}

In contrast, our main theorem 
admits ramification.
We just assume that 
the ramification index $e$ of a fixed splitting field $L/F$
is relatively prime to
the degree $n$ of the cover.
Precisely:

\begin{mth}[\ref{tMt}]
  Let $F$ be a non-archimedean local field of characteristic zero,
  $\mathbf f$ the residue field of the valuation ring of $F$, and
  $\mu_n\subset \mathbf f^\times$ the cyclic subgroup of order $n$.
  Suppose that an algebraic torus $T$ defined over $F$ 
  splits over a finite Galois extension $L/F$ 
  whose ramification index $e$ is relatively prime to $n$.
  Let $\mu_n\to\widetilde{T}\to T$ be a Brylinski-Deligne 
  covering group.

  Then the group $\mathcal S=\mathcal S_{\widetilde T}$ 
  parametrizing a packet of representations of $\widetilde{T}$ is 
  written as the quotient
  \[
    \mathcal S_{\widetilde T}
    =\left.\iota((Y^{\Gamma\#}\otimes_{\mathbb Z}
    \overline{\mathbf f}^\times)^\Gamma)\middle/
    \iota((Y^{\#}\otimes_{\mathbb Z}
    \overline{\mathbf f}^\times)^\Gamma)\right..
  \]
  of subgroups of 
  $(Y\otimes_{\mathbb Z}\overline{\mathbf f}^\times)^\Gamma$.
\end{mth}

As already mentioned,
the main theorem generalizes
Weissman's work on covering groups.
Moreover, this theorem partially realizes his 
hope~\cite[Remark 5.18]{Weiss_loc_tori} to parametrize a packet
in the 
correspondence for a cover of a ramified torus.

The next two sections are devoted to formulating our subjects,
namely a Brylinski-Deligne covering group $\widetilde T\to T$ of
a torus
and its local Langlands correspondence.
Especially, we define the finite group 
$\mathcal S_{\widetilde T}$ as a quotient of the image
$\mathrm Z^\dagger\subset T$ of the center 
$\operatorname Z(\widetilde T)\subset\widetilde T$
in \Cref{eLlc}.
To prove the main theorem, we describe 
the subgroup $\mathrm Z^\dagger\subset T$ 
by applying the Galois cohomology theory in \Cref{sZd}.
In \Cref{Weiss},
we overcome a critical part by reducing it to a certain orthogonality
in the local Tate duality.
We then prove the main theorem in the final section.

\subsection*{Acknowledgements}
The author would like to thank his advisor Professor 
Tamotsu Ikeda for fruitful discussions, 
and thank Miyu Suzuki for providing me with valuable information.
This work was supported by JSPS KAKENHI Grant Number JP23KJ1298.

\subsection*{Notation}
Let $\mathbb Z_{>0}$, $\mathbb Z$, $\mathbb Q$, and $\mathbb C$
respectively denote
the set of positive rational integers,
the ring of rational integers,
the field of rational numbers, and
the field of complex numbers.
The symbol $\otimes$ always denotes the tensor product 
over $\mathbb Z$.
For a map $f\colon A\to B$ and a subset $A'\subset A$,
let $f|_{A'}$ denote the restriction, as usual.
For a group $G$, let $\operatorname Z(G)\subset G$ denote the center.
If a group $G$ acts on an abelian group $A$, then
$A^G=\operatorname H^0(G,A)$ is the subgroup consisting of
fixed points.
For a ring $R$, let $R^\times$ denote the group of units.

We write $F$ for a non-archimedean local field
of characteristic zero,
$\mathcal O=\mathcal O_F$ for its valuation ring, and
$\mathbf f$ for the residue field.
We fix an algebraic closure $\overline F/F$, and write
$\overline{\mathcal O}\subset\overline F$ for its valuation ring
and $\overline{\mathbf f}$ for the residue field of
$\overline{\mathcal O}$.
Note that $\overline{\mathbf f}$ is an algebraic closure of 
the field $\mathbf f$.
For a field $k$, let 
$\Gamma_k=\operatorname{Gal}(\overline k/k)$ denote
its absolute Galois group.
In the case of the local field $F$, we often omit the subscript
and simply write $\Gamma=\Gamma_F$.

In this paper, we treat an algebraic torus $\mathbb T$
defined over the local field $F$, and 
fix a finite Galois extension $L/F$ where $\mathbb T$ splits over $L$.
Let $Y=\operatorname{Hom}(\mathbb G_{\mathrm m},\mathbb T)$ be
the cocharacter lattice, 
$T=\mathbb T(F)=(Y\otimes_{\mathbb Z}\overline F^\times)^\Gamma$ 
the group of rational points, and 
$e\in\mathbb Z_{>0}$ the ramification index of the extension $L/F$.
We write $\operatorname{ord}\colon F^\times\to\mathbb Z$ for
the map of additive valuation, and extend it to
$\overline F^\times\to\mathbb Q$.
We continue to write $\operatorname{ord}$ for the induced map
\[
  T=(Y\otimes\overline F^\times)^\Gamma\to(Y\otimes\mathbb Q)^\Gamma
  =Y^\Gamma\otimes\mathbb Q.
\]
For a subgroup $Y'\subset Y$, let $\iota$ denote
the inclusion map $Y'\hookrightarrow Y$
and the induced ones.
For instance,
$\iota\colon(Y'\otimes\overline F^\times)^\Gamma\to
(Y\otimes\overline F^\times)^\Gamma$
is a morphism between tori.

\section{Brylinski-Deligne covering groups of a torus}
\label{eCov} 
In the literature on Langlands correspondence for
general covering groups~\cite{GG_Weiss,Weiss_L}, 
the authors often focus on certain covering groups defined by
Brylinski and Deligne~\cite{Bry_Del}.
In this section, we review the definition and 
the classification theorem for Brylinski-Deligne covering groups
of an algebraic torus $\mathbb T$ defined over a local field $F$.

Roughly speaking, a Brylinski-Deligne covering group is
a central extension by Milnor's group $\operatorname K_2$
in algebraic K-theory, or its variation.
The groups $\operatorname{K}_2$ are 
generalized by Quillen~\cite{Qill_hi}
to be defined over schemes, and form a presheaf on the big Zariski
site of the base field $F$.
Let $\mathbb K_2$ be the Zariski sheaf associated with the presheaf 
$\operatorname K_2$.

To define a Brylinski-Deligne covering group,
let $\mu_n
\subset F^\times$ denote the cyclic subgroup of
order $n$.

\begin{dfn}[Brylinski-Deligne~\cite{Bry_Del}]
  \label{d_BD}
  Let $\mathbb T\to\operatorname{Spec}F$ be an algebraic torus.
  According to the context, a \emph{Brylinski-Deligne covering group}
  of $\mathbb T$ or $T=\mathbb T(F)$ means one of the following:
  \begin{enumerate}
    \item\label{e_BD_sh} a central extension
      $1\to\mathbb{K}_2\to\widetilde{\mathbb{T}}\to\mathbb{T}\to 1$
      as a sheaf of groups on the big Zariski site of 
      $\operatorname{Spec}F$,
    \item\label{e_BD_sect} the section 
      $1\to\operatorname K_2(F)\to\widetilde{\mathbb{T}}(F)\to{T}
      \to 1$ of~\ref{e_BD_sh} on $\operatorname{Spec}F$, or
    \item\label{e_BD_topl} the push-out 
      $1\to\mu_n\to\widetilde{T}\to{T}\to 1$ of~\ref{e_BD_sect}
      via the Hilbert symbol
      $\operatorname K_2(F)\to\mu_n$. 
  \end{enumerate}
\end{dfn}

A Brylinski-Deligne covering group defined 
in \ref{d_BD}~\ref{e_BD_topl} is indeed a topological covering
group~\cite[Construction 10.3]{Bry_Del}.
Further, Brylinski and Deligne generally defined their covering groups
for any algebraic group.

For a split torus $T$, a 2-cocycle of a Brylinski-Deligne cover
$\mu_n\to\widetilde{T}\to{T}$ is written as a product of $n$-th
Hilbert symbols 
$(\,,)_n\colon F^\times\times F^\times\to\mu_n$~\cite[3.3]{GG_Weiss}.
For instance,
let $T=F^\times$ be the one-dimensional split torus, i.e.
just the multiplicative group.
Then for an integer $a$, the 2-cocycle 
$(\,,)_n^a\colon F^\times\times F^\times\to\mu_n$ determines
a Brylinski-Deligne cover
$\mu_n\to\widetilde{T}\to{T}$.
As another example,
let $T=F^\times\times F^\times$ be the two-dimensional
split torus.
Then for four integers $(a_{ij})_{i,j\in\{1,2\}}$, the 2-cocycle
\[
  (F^\times \times F^\times)\times(F^\times \times F^\times)\to
  \mu_n,\ ((s_1,s_2),(t_1,t_2))\mapsto\prod_{i,j\in\{1,2\}}
  (s_i,t_j)_n^{a_{ij}}
\]
gives a Brylinski-Deligne cover
$\mu_n\to\widetilde{T}\to{T}$.
Especially if $n\geq 2$ and 
$(a_{ij})_{i,j}=
\begin{pmatrix}
  0&1\\
  0&0
\end{pmatrix}
$,
then the covering group $\widetilde T$ is non-abelian. 

For any torus $T$, one has a finite Galois extension $L/F$
of the base field such that $T$ splits over $L$.
Then the following classification theorem gives a bilinear form
$Y\times Y\to\mathbb Z$ on the cocharacter lattice,
which is a key ingredient to describe our main theorem.

\begin{thm}[Classification Theorem of covers~\cite{Bry_Del}]
  \label{etBil}
  The following two Picard categories are equivalent:
  \begin{enumerate}
    \item the category of Brylinski-Deligne covers
      $\mathbb{K}_2\to\widetilde{\mathbb{T}}\to\mathbb{T}$, and 

    \item the category of the following data $(Q,\mathcal{E})$:
      \begin{itemize}
        \item 
          $Q\colon Y\to\mathbb{Z}$ is a Galois-invariant quadratic 
          form on the cocharacter lattice, 
          which gives a bilinear form 
          $B\colon Y\times Y\to\mathbb{Z}$ by 
          \[(y,y')\mapsto Q(y+y')-Q(y)-Q(y').\]
        \item 
          $1\to L^\times\to\mathcal{E}\stackrel{\overline{(\ )}}{\to}Y
          \to 1$ is a Galois-equivariant {central extension} such that
          for any $\epsilon,\zeta\in\mathcal{E}$,
          \[\epsilon\zeta\epsilon^{-1}\zeta^{-1}=(-1)^{B(\overline{\epsilon},\overline{\zeta})}\]
          in $L^\times$.
        \item
          For any data $(Q,\mathcal E)$ and $(Q',\mathcal E')$,
          \[
            \operatorname{Hom}((Q,\mathcal{E}),(Q',\mathcal{E}'))=
            \begin{cases}
              \emptyset &\text{if }Q\neq Q'\\
              \{\text{morphisms }\mathcal{E}\to\mathcal{E}'
              \text{ of extensions}\}
              &\text{if }Q=Q'.
            \end{cases}
          \]
      \end{itemize}
  \end{enumerate}
\end{thm}

Brylinski and Deligne generally proved this type of classification
theorem for covering groups of reductive algebraic 
groups~\cite{Bry_Del}.

\begin{dfn}
  \label{eBil}
  Let $\mathbb K_2\to\widetilde{\mathbb T}\to\mathbb T$ be
  a Brylinski-Deligne cover, and
  $Q_{\widetilde{\mathbb T}}\colon Y\to\mathbb Z$
  the Galois-invariant quadratic form corresponding to 
  $\widetilde{\mathbb T}$ by \Cref{etBil}.
  Then
  we write $B_{\widetilde{\mathbb T}}\colon Y\times Y\to\mathbb{Z}$ 
  for the associated bilinear form 
  $(y,y')\mapsto Q(y+y')-Q(y)-Q(y')$.
  If $\mu_n\to\widetilde T\to T$ is the cover induced from
  $\widetilde{\mathbb T}$ (\Cref{d_BD}), then we write
  $
  B_{\widetilde T}=B_{\widetilde{\mathbb T}}$ by abuse of notation.
\end{dfn}

\section{Local Langlands correspondence for a covering group}
\label{eLlc}
In this section, we introduce a finite abelian group 
$\mathcal S=\mathcal S_{\widetilde T}$, which is the main object of this work
(\Cref{Das}).
Then we shall explain how $\mathcal S_{\widetilde T}$ parametrizes a packet in 
the local Langlands correspondence for a Brylinski-Deligne 
covering group $\mu_n\to\widetilde T\to T$ of a 
torus.
Our formulation is essentially due to Weissman~\cite{Weiss_loc_tori}.

The local Langlands correspondence for a covering group 
$\widetilde T\to T$ of a torus aims to parametrize 
genuine irreducible representations of $\widetilde T$.
To treat genuine representaions, we throughout fix an embedding 
\[j\colon\mu_n\to\mathbb{C}^\times\]
of the cyclic group of order $n$.
Recall that
a representation of a covering group $\mu_n\to\widetilde G\to G$ 
of a topological group is \emph{genuine} 
if $\mu_n\subset\widetilde{G}$ acts as 
multiplication by scalers in $\mathbb{C}$ via the embedding $j$.
That is, for any $\zeta\in\mu_n$ and any vector 
$v$ in the representation space, one has
\[\zeta\cdot v=j(\zeta)v.\]
A \emph{genuine character} 
$\chi\colon\widetilde G\to\mathbb C^\times$ is
a character that is genuine as a representation,
i.e. $\chi|_{\mu_n}=j$.

\begin{exm}[central character]
  \label{xCc}
  \begin{enumerate}
    \item
      Let $\mu_n\to\widetilde T\to T$ be a covering group of a torus,
      and $\rho$ an irreducible genuine representation of 
      $\widetilde T$.
      Let $\mathrm Z^\dagger\subset T$ denote the image of the center 
      $\operatorname Z(\widetilde T)\subset \widetilde T$, 
      and note that
      \[\mu_n\to\operatorname Z(\widetilde T)\to\mathrm Z^\dagger\]
      is a covering group.
      Then the central character 
      $\rho|_{\operatorname Z(\widetilde T)}\colon
      \operatorname Z(\widetilde T)\to\mathbb C^\times$
      is a genuine character.
    \item\label{eSvN}
      Let $\mu_n\to\widetilde T\to T$ be a Brylinski-Deligne covering
      of a torus.
      Then the correspondence from the genuine irreducible
      representations of the covering group $\widetilde T$ to
      their central characters 
      is bijective onto the set of genuine characters 
      $\operatorname Z(\widetilde T)\to
      \mathbb C^\times$~\cite[Theorem 3.1]{Weiss_loc_tori}.
  \end{enumerate}
\end{exm}


To establish the local Langlands correspondence for 
the Brylinski-Deligne covering group $\widetilde T\to T$, 
it suffices to parametrize the genuine characters on the center
$\operatorname Z(\widetilde T)$, 
by \Cref{xCc}~\ref{eSvN}.
To estimate $\operatorname Z(\widetilde T)\subset\widetilde T$,
we choose the following isogeny $T^\#\stackrel\iota\to T$ of tori,
which approximates the image $\mathrm Z^\dagger\subset T$ of
$\operatorname Z(\widetilde T)$.

\begin{dfn}
  \label{tTs}
  Let $Y=\operatorname{Hom}(\mathbb G_{\mathrm m},\mathbb T)$ be 
  the cocharacter group of the torus $T=\mathbb T(F)$.
  Let $B_{\widetilde T}\colon Y\times Y\to\mathbb Z$
  be the bilinear form attached to the covering 
  $\mu_n\to\widetilde T\to T$, precisely 
  to the Brylinski-Deligne cover 
  $\mathbb K_2\to\widetilde{\mathbb T}\to\mathbb T$ defining
  $\widetilde T$ (\Cref{eBil}).
  \begin{enumerate}
    \item
      We define the subgroup $Y^\#\subset Y$ by
      \[Y^\#=\{y\in Y\mid B_{\widetilde T}(y,Y)\subset n\mathbb Z\}.\]
    \item
      We write $T^\#=(Y^\#\otimes\overline{F}^\times)^\Gamma$ 
      for the torus defined by $Y^\#$.
      Then the inclusion $Y^\#\hookrightarrow Y$ induces an isogeny
      $\iota\colon T^\#\to T$.
  \end{enumerate}
\end{dfn}

It is known that $\iota(T^\#)\subset\mathrm{Z}^\dagger$, 
and that this index is finite~\cite[Theorem 1.3]{Weiss_LG_tori}.

\begin{dfn}
  \label{Das}
  For the Brylinski-Deligne covering group $\widetilde T\to T$,
  we define the finite abelian group 
  $\mathcal S=\mathcal S_{\widetilde T}$ as the quotient
  \[\mathcal S_{\widetilde T}=\mathrm{Z}^\dagger/\iota(T^\#).\]
  This may be called the \emph{packet group}~\cite{Weiss_loc_tori}.
\end{dfn}

In the rest of this section, we explain the role of this group
$\mathcal S_{\widetilde T}$.
First, $\mathcal S_{\widetilde T}$ and the other groups defined above
form the 
exact sequences
\[
  \begin{array}{*9c}
    && 1\\
    && \uparrow\\
    && \mathcal S_{\widetilde T}\\
    && \uparrow\\
    1 &\gets& \mathrm Z^\dagger &\gets& 
    \operatorname Z(\widetilde T) &\gets& \mu_n &\gets& 1
    \\
    && 
    \uparrow
    \\
    1 &\gets& \iota(T^\#) &\stackrel\iota\twoheadleftarrow& T^\#
    \\
    && \uparrow\\
    && 1\rlap.
  \end{array}
\]
To apply the functor $\operatorname{Hom}(\quad,\mathbb{C}^\times)$
to this diagram, we note the following facts.

\begin{fct}
  \begin{enumerate}
    \item
      The multiplicative group $\mathbb C^\times$ is 
      an injective object in the category of locally compact abelian 
      groups~\cite[Th\'eor\`eme 5]{Dix_lca}.
      That is, for any locally compact abelian group $A$ and 
      its closed subgroup $C\subset A$, any continuous character
      $C\to\mathbb C^\times$ extends to a continuous character
      $A\to\mathbb C^\times$.
    \item
      Let $\operatorname{W}_F$ be the Weil group of 
      the local field $F$,
      and $\widehat{T^\#}=\operatorname{Hom}(Y^\#,\mathbb C^\times)$
      the Langlands dual of the torus $T^\#$.
      Then 
      the Langlands correspondence for 
      $T^\#$~\cite[Th\'eor\`eme 6.2]{Lab_Lgp} asserts an isomorphism
      \[
        \operatorname{LLC}_{T^\#}\colon
        \operatorname{Hom}(T^\#,\mathbb{C}^\times)\cong
        \operatorname{H}^1(\operatorname{W}_F,\widehat{T^\#}).
      \]
  \end{enumerate}
\end{fct}

These facts ensures the following frame of the local 
Langlands correspondence for the covering group $\widetilde T\to T$.

\begin{prp}
  \label{tLlc}
  The objects introduced above fit in the diagram
  \[
    \begin{array}{*9c}
      1 && 
      \hbox to0pt{\hss$
      \{\text{genuine irreducible representations of }\widetilde T\}
      $\hss}
      \\
      \downarrow && 
      \llap{$\scriptstyle (\;)|_{\operatorname Z(\widetilde T)}$}
      \downarrow
      \rlap{$\scriptstyle 1:1$}
      \\
      \operatorname{Hom}(\mathcal S_{\widetilde T},\mathbb{C}^\times)
      && \{\text{genuine characters}\} &\to& \{j\}
      \\
      \downarrow && \cap && \cap
      \\
      \llap{$1\to$\,}
      \operatorname{Hom}(\mathrm{Z}^\dagger,\mathbb{C}^\times)
      &\to&
      \operatorname{Hom}(\operatorname{Z}(\widetilde T),
      \mathbb{C}^\times)
      &\to&
      \operatorname{Hom}(\mu_n,\mathbb{C}^\times) 
      \rlap{\;$\to1$}
      \\
      \mbox{}
      \downarrow
      \rlap{$\scriptstyle(\;)|_{\,\iota(T^\#)}$}
      \\
      \operatorname{Hom}(\iota(T^\#),\mathbb{C}^\times)
      &\stackrel{\iota^*}\hookrightarrow&
      \operatorname{Hom}(T^\#,\mathbb{C}^\times)
      \\
      \downarrow &&
      \rB{-90}{$\cong$}
      \rlap{\,$\scriptstyle\operatorname{LLC}_{T^\#}$}
      \\
      1 &&
      \operatorname{H}^1(\operatorname{W}_F,\widehat{T^\#}),
    \end{array}
  \]
  where the middle horizontal sequence and the left vertical one
  are exact.
\end{prp}

By definition, the group
$\operatorname{Hom}(\mathrm{Z}^\dagger,\mathbb{C}^\times)$
acts simply transitive on the set of genuine characters
$\operatorname Z(\widetilde T)\to\mathbb C^\times$.
Thus, 
fixing a genuine character $\chi$ on 
$\operatorname Z(\widetilde T)$ as a base point defines
a bijection
\[
  b_\chi\colon
  \{\text{genuine characters on }\operatorname Z(\widetilde T)\}
  \to
  \operatorname{Hom}(\mathrm{Z}^\dagger,\mathbb{C}^\times).
\]
By \Cref{tLlc}, the composed map
\[
  \operatorname{LLC}_{T^\#}\circ\iota^*\circ(\,)|_{\iota(T^\#)}
  \colon\operatorname{Hom}(\mathrm{Z}^\dagger,\mathbb{C}^\times)
  \to\operatorname{H}^1(\operatorname{W}_F,\widehat{T^\#})
\]
has kernel
$\operatorname{Hom}(\mathcal S_{\widetilde T},\mathbb{C}^\times)$,
which is a finite 
group.

\begin{dfn}
  \label{iLc}
  Fix a genuine character 
  $\chi\colon\operatorname Z(\widetilde T)\to\mathbb C$.
  Then the \emph{local Langlands correspondence} for 
  the Brylinski-Deligne cover 
  $\widetilde T\to T$ is the composition
  \[
    \phantom{\text{genuine}}
    \begin{array}[b]{*9c}
      \hbox to0pt{\hss$
      \{\text{genuine irreducible representations of }
      \widetilde T\}
      $\hss}\\
      \llap{$\scriptstyle b_\chi\circ
      (\,)|_{\operatorname Z(\widetilde T)}$}
      \downarrow\rlap{$\scriptstyle 1:1$}\\
      \operatorname{Hom}(\mathrm{Z}^\dagger,\mathbb{C}^\times)
    \end{array}
    \xrightarrow{\operatorname{LLC}_{T^\#}\circ\iota^*\circ
    (\,)|_{\iota(T^\#)}}
    \operatorname{H}^1(\operatorname{W}_F,\widehat{T^\#}).
  \]
\end{dfn}

As remarked above, this is a finite-to-one correspondence.
Note that some packets, i.e., fibers, may be empty.
By definition, the non-empty packets are described by the group
$\mathcal S_{\widetilde T}$, as follows.

\begin{prp}[role of $\mathcal S_{\widetilde T}$]
    \label{tlS}
      In the local Langlands correspondence defined in~\ref{iLc}, 
      each non-empty packet $\Pi$ is just an orbit under the group 
      $\operatorname{Hom}(\mathcal S_{\widetilde T},\mathbb C^\times)$.
      Thus,
      Choosing a base point in the packet $\Pi$ defines
      a one-to-one correspondence between $\Pi$ and 
      $\operatorname{Hom}(\mathcal S_{\widetilde T},\mathbb C^\times)$.
\end{prp}
The correspondence between a packet $\Pi$ and the group 
$\operatorname{Hom}(\mathcal S_{\widetilde T},\mathbb C^\times)$
stated in \Cref{tlS} is an analogue of 
the usual formulation to parametrize a packet in the local Langlands 
correspondence for a reductive algebraic group.

\section{Description of the image $\mathrm{Z}^\dagger$ of the center}
\label{sZd}
To determine the finite group 
$\mathcal S_{\widetilde T}=\mathrm Z^\dagger/\iota(T^\#)$
parametrizing a packet,
we prepare exact sequences, 
and describe the subgroup $\mathrm{Z}^\dagger\subset T$
via a non-degenerate bilinear form.

For $k\in\mathbb Z$ and
a finite 
$\Gamma$-module $M$,
let $M(k)$ denote the $k$-th Tate twist of $M$.
For instance,
$(\mathbb Z/m\mathbb Z)(1)=\mu_m\subset\overline F^\times$
is the cyclic subgroup of order $m$, and
$(\mathbb Z/m\mathbb Z)(2)=\mu_m\otimes\mu_m$,
for $m\in\mathbb Z_{>0}$.

\begin{lmm}
  \label{lYY}
  Let $A$ be a $\Gamma$-module, and
  $Y'\subset Y$ a $\Gamma$-submodule of the cocharacter lattice.
  Assume that for some $m\in\mathbb Z_{>0}$, we have the inclusion
  $mY\subset Y'$ 
  and exact sequence
  \[0\to(\mathbb{Z}/m\mathbb{Z})(1)\to A\stackrel{m}{\to}A\to 0.\]
  Then they induce a short exact sequence
  \[0\to(Y/Y')(1)\to Y'\otimes A\to Y\otimes A\to 0.\]
\end{lmm}

\begin{proof}
  The sequence $0\to Y'\to Y\to Y/Y'\to 0$ induces an exact sequence
  \[
    0\to \operatorname{Tor}_1(Y/Y',A)\to 
    Y'\otimes A\to Y\otimes A \to 0.
  \]
  It suffices to give an isomorphism
  $\operatorname{Tor}_1(Y/Y',A)\cong(Y/Y')(1)$. 
  Indeed, the sequence
  $0\to(\mathbb{Z}/m\mathbb{Z})(1)\to A\stackrel{m}{\to} A\to 0$ 
  gives an exact sequence
  \[
    \begin{matrix}
      \operatorname{Tor}_{1}(Y/Y',A)\stackrel{0}{\to}
      \operatorname{Tor}_{1}(Y/Y',A)&\to(Y/Y')\otimes(\mathbb{Z}/m\mathbb{Z})(1)\to&(Y/Y')\otimes A\\
      &\rB{90}{=}                                    &\rB{90}{=}\\
      &(Y/Y')(1)                                     &0.
    \end{matrix}
  \]
\end{proof}

\begin{exm}
  \label{gSf}
  \begin{enumerate}
    \item\label{eSf}
      Let $A=\overline{\mathbf f}^\times$ be an algebraic closure
      of the residue field, and 
      $Y'=Y^\#=\{y\in Y\mid B_{\widetilde T}(y,Y)\in n\mathbb Z\}$
      as defined in \ref{tTs}.
      Suppose the inclusion 
      $\mu_n\subset\mathbf f^\times$, so that
      the assumption in \Cref{lYY} holds for $m=n$.
      Then we have a short exact sequence
      \[
        0\to(Y/Y^\#)(1)\to Y^\#\otimes \overline{\mathbf f}^\times\to
        Y\otimes \overline{\mathbf f}^\times\to 0.
      \]

    \item \label{eOF}
      Let $Y'\subset Y$ be a $\Gamma$-submodule of the same rank,
      i.e.
      some $m\in\mathbb Z_{>0}$ satisfies $mY\subset Y'$.
      Then the pairs $(\overline{\mathcal O}^\times,Y')$ and
      $(\overline F^\times,Y')$ respectively satisfy the assumption
      in \Cref{lYY}. 
      Hence the tensor products with the valuation sequence
      $1\to\overline{\mathcal{O}}^\times\to\overline{F}^\times
      \stackrel{\operatorname{ord}}{\to}\mathbb{Q}\to 0$
      gives the commutative diagram
      \[
        \begin{array}{ccccc}
          &0&&0&\\
          &\rB{-90}{$\to$}&&\rB{-90}{$\to$}
          &\\
          0\to(Y/Y')(1)\to&Y'\otimes
          \overline{\mathcal{O}}^\times&
          \to&Y\otimes\overline{\mathcal{O}}^\times&\to 0\\
          \rB{90}{$=$}&\rB{-90}{$\to$}&&\rB{-90}{$\to$}
          &\\
          0\to(Y/Y')(1)\to&Y'\otimes\overline{F}^\times&\to&
          Y\otimes\overline{F}^\times&\to 0\\
          &\rB{-90}{$\to$}\rlap{\,$\scriptstyle\operatorname{ord}$}
          &&\rB{-90}{$\to$}\rlap{\,$\scriptstyle\operatorname{ord}$}
          &\\
          &Y'\otimes\mathbb{Q}&\stackrel{\sim}{\to}& 
          Y\otimes\mathbb{Q}\\
          &\rB{-90}{$\to$}&&\rB{-90}{$\to$}
          &\\
          &0&&0&
        \end{array}
      \]
      of exact sequences.
  \end{enumerate}
\end{exm}

In \Cref{gSf}~\ref{eOF},
note that $(Y'\otimes\mathbb Q)^\Gamma=Y'^\Gamma\otimes\mathbb Q$.
Then, the Galois cohomology sequences are described as follows.

\begin{prp}
  \label{ofq_Gal}
  Let $Y'\subset Y$ be a $\Gamma$-submodule of the same rank.
  Then the tensor products with the valuation sequence
  $1\to\overline{\mathcal{O}}^\times\to\overline{F}^\times
  \stackrel{\operatorname{ord}}{\to}\mathbb{Q}\to 0$
  gives the 
  diagram 
  \[
    \begin{matrix}
      &0                                                   &                     &0                                                                                                              \\
      &\rB{-90}{$\to$}                                       &                     &\rB{-90}{$\to$}                                                                                                  \\
      0\to(Y/Y')(1)^\Gamma\to&(Y'\otimes\overline{\mathcal{O}}^\times)^\Gamma   &\stackrel{\iota}{\to}&(Y\otimes\overline{\mathcal{O}}^\times)^\Gamma      &\stackrel{\partial_{Y/Y'}}{\to} \operatorname{H}^1(\Gamma,(Y/Y')(1))\\
      \rB{90}{$=$}&\rB{-90}{$\to$} &                                                    &\rB{-90}{$\to$}        &\rB{90}{$=$}                                          \\
      0\to(Y/Y')(1)^\Gamma\to&(Y'\otimes\overline{F}^\times)^\Gamma             &\stackrel{\iota}{\to}&(Y\otimes\overline{F}^\times)^\Gamma                &\stackrel{\partial_{Y/Y'}}{\to}\operatorname{H}^1(\Gamma,(Y/Y')(1))\\
      &\rB{-90}{$\to$}\rlap{\,$\scriptstyle\operatorname{ord}$}&                     &\rB{-90}{$\to$}\rlap{\,$\scriptstyle\operatorname{ord}$}                                                           \\
      &Y'^{\Gamma}\otimes\mathbb{Q}                       &\stackrel{\sim}{\to} &Y^\Gamma\otimes\mathbb{Q}.                          &
    \end{matrix}
  \]
  of exact sequences.
\end{prp}

We write $\partial=\partial_{Y/Y'}$ for the connecting homomorphisms.
Let
$V_{Y'}=\operatorname{ord}((Y'\otimes\overline{F}^\times)^\Gamma)$
be the image of $\operatorname{ord}$ in the diagram of
\Cref{ofq_Gal}.

\begin{crr}
  \label{ePbi}
  \begin{enumerate}
    \item\label{eCdi}
      The images of $\partial$ and $\operatorname{ord}$ in
      \Cref{ofq_Gal} gives the commutative diagram 
      \[
        \begin{matrix}
          &                  &   &0                                                &                     &0                                             &   &0                                                       &\\
          &                  &   &\rB{-90}{$\to$}                                    &                     &\rB{-90}{$\to$}                                 &   &\rB{-90}{$\to$}                                           &\\
          0\to&(Y/Y')(1)^\Gamma&\to&(Y'\otimes\overline{\mathcal{O}}^\times)^\Gamma&\stackrel{\iota}{\to}&(Y\otimes\overline{\mathcal{O}}^\times)^\Gamma&\stackrel\partial\to&\partial((Y\otimes\overline{\mathcal{O}}^\times)^\Gamma)&\to 0\\
          &\rB{-90}{$=$}       &   &\rB{-90}{$\to$}                                    &                     &\rB{-90}{$\to$}                                 &   &\rB{-90}{$\hookrightarrow$}                                           &\\
          0\to&(Y/Y')(1)^\Gamma&\to&(Y'\otimes\overline{F}^\times)^\Gamma          &\stackrel{\iota}{\to}&(Y\otimes\overline{F}^\times)^\Gamma          &\stackrel\partial\to&\partial((Y\otimes\overline{F}^\times)^\Gamma)          &\to 0\\
          &                  &   &\rB{-90}{$\to$}                                    &                     &\rB{-90}{$\to$}\rlap{\,$\scriptstyle\operatorname{ord}$}                                 &   &\rB{-90}{$\to$}                                           &\\
          &0                 &\to&V_{Y'}                                         &\hookrightarrow                  &V_Y                                           &\to&V_Y/V_{Y'}                                            &\to 0\\
          &                  &   &\rB{-90}{$\to$}                                    &                     &\rB{-90}{$\to$}                                 &   &\rB{-90}{$\to$}                                           &\\
          &                  &   &0                                                &                     &0                                             &   &0                                                      &
        \end{matrix}
      \]
      of exact sequences.

    \item\label{eePbi}
      In the diagram~\ref{eCdi}, we have
      \[
        \iota((Y'\otimes\overline{F}^\times)^\Gamma)\cap
        (Y\otimes\overline{\mathcal{O}}^\times)^\Gamma
        =\iota((Y'\otimes\overline{\mathcal{O}}^\times)^\Gamma).
      \]
    \item \label{eDsj}
      Suppose that $T'=(Y'\otimes\overline F^\times)^\Gamma$ is 
      a split torus over $F$, i.e. $\Gamma$ acts trivially on $Y'$.
      Then
      $\partial
      ((Y\otimes\overline F^\times)^\Gamma)
      =\operatorname{H}^1(\Gamma,(Y/Y')(1))$.
  \end{enumerate}
\end{crr}

\begin{proof}
  \ref{eCdi}
  This follows from the snake lemma.
  \ref{eePbi}
  The diagram displays that the map
  \[
    \left.(Y\otimes\overline{\mathcal O}^\times)^\Gamma\middle/
    \iota((Y'\otimes\overline{\mathcal O}^\times)^\Gamma)\right.\to
    \left.(Y\otimes\overline{F}^\times)^\Gamma\middle/
    \iota((Y'\otimes\overline{F}^\times)^\Gamma)\right.
  \]
  is injective.
  \ref{eDsj}
  Hilbert's theorem 90 shows 
  $\operatorname{H}^1(\Gamma,Y'\otimes\overline{F}^\times)=0$.
\end{proof}

For the subgroup $Y'=nY$ of $Y$, the connecting homomorphism
\[
  \partial=\partial_{Y/nY}\colon
  T=(Y\otimes\overline F^\times)^\Gamma\to
  \operatorname{H}^1(\Gamma,(Y/nY)(1))
\]
is applied to express commutators of lifts of 
elements in $T$ to $\widetilde T$.

\begin{fct}[{\cite[Proposition 9.9]{Bry_Del}}]
  \label{tCmn}
  The composite map 
  \[
    T\times T\stackrel{\partial\times\partial}{\to}
    \operatorname{H}^1(\Gamma,(Y/nY)(1))^{\oplus2}
    \stackrel{\cup}{\to}\operatorname{H}^2(\Gamma,
    (Y/nY)(1)^{\otimes2})\stackrel{B_{\widetilde T}}{\to}
    \operatorname{H}^2(\Gamma,\mathbb{Z}/n(2))=\mu_n
  \]
  gives the commutator of 
  lifts of two elements in $T$ to $\widetilde T$.
  That is, for any $s,t\in T$ and their lift 
  $\tilde s,\tilde t\in\widetilde T$, we have
  \[
    \operatorname{comm}_{\widetilde T}(s,t):=
    B_{\widetilde T}(\partial s\cup\partial t)=
    \tilde s\tilde t\tilde s^{-1}\tilde t^{-1}.
  \]
  Thus,
  $\mathrm{Z}^\dagger=\{t\in T\mid
  \forall t'\in T,B_{\widetilde T}(\partial t\cup\partial t')=0\}$.
\end{fct}

In the rest of this section,
we slightly generalize the description of the commutator map
$\operatorname{comm}_{\widetilde T}$ in \Cref{tCmn},
so that we can apply other useful connecting homomorphisms.

For $i=1,2$, 
let $T_i$ be an algebraic torus with cocharacter lattice $Y_i$.
\begin{dfn}
  \label{dMT}
  Let $B\colon Y_1\times Y_2\to\mathbb Z$ be a $\Gamma$-invariant 
  bilinear form.
  Then
  we define the $\Gamma$-submodule 
  $Y_2^\#\subset
  Y_1$ as the annihilator
  \[Y_2^\#=\{y\in Y_1\mid B(y,Y_2)\subset n\mathbb Z\}\]
  of $Y_2$ for the bilinear form $Y_1\times Y_2\to\mathbb Z/n$.
  We similarly define the $\Gamma$-submodule 
  $Y_1^\#\subset Y_2$ by
  $Y_1^\#=\{y\in Y_2\mid B(Y_1,y)\subset n\mathbb Z\}$.
\end{dfn}

For $B=B_{\widetilde T}\colon Y\times Y\to\mathbb Z$,
\Cref{dMT} agrees with the definition of the submodule $Y^\#$ 
in \ref{tTs}.
In general, the modules $Y_2^\#$ and $Y_1^\#$ are used 
to define non-degenerate bilinear forms:

\begin{prp}
  \label{eNdg}
  Let $B\colon Y_1\times Y_2\to\mathbb Z$ be
  a $\Gamma$-invariant bilinear form.
  For $i=1,2$, let $Y_i'\subset Y_i$ be a submodule.
  \begin{enumerate}
    \item\label{eQb}
      The form $B$ induces a bilinear form
      $Y_1/Y_1'\times Y_2/Y_2'\to\mathbb Z/n$ if and only if
      $Y_1'\subset Y_2^\#$ and $Y_2'\subset Y_1^\#$.
    \item\label{eNdm}
      The induced bilinear form
      $B\colon Y_1/Y_2^\#\times Y_2/Y_1^\#\to\mathbb Z/n$ 
      for the subgroups $Y_2^\#$ and $Y_1^\#$ is non-degenerate.
    \item\label{eTd}
      The bilinear form 
      $B\colon Y_1/Y_2^\#\times Y_2/Y_1^\#\to\mathbb Z/n$
      induces a non-degenerate bilinear form 
      \[
        \begin{array}{*9c}
          \operatorname{H}^1(\Gamma,(Y_1/Y_2^\#)(1))\times
          \operatorname{H}^1(\Gamma,(Y_2/Y_1^\#)(1))
          &\to&\operatorname{H}^2(\Gamma,\mathbb{Z}/n(2))\\
          \rB{90}{$\in$} && \rB{90}{$\in$}\\
          (h_1,h_2) &\mapsto& B(h_1\cup h_2).
        \end{array}
      \]
  \end{enumerate}
\end{prp}

\begin{proof}
  \ref{eQb} and \ref{eNdm}
  hold by definition.

  \ref{eTd}
  By definition, the bilinear form
  $B\colon(Y_1/Y_2^\#)(1)\times(Y_2/Y_1^\#)(1)\to\mathbb{Z}/n(2)\cong
  \mathbb{Z}/n(1)$
  is non-degenerate, i.e., gives a pairing of dual Galois modules.
  By local Tate duality~\cite[II.5.2]{Srr_Gal_coh}, the cup product
  $\operatorname{H}^1(\Gamma,(Y_1/Y_2^\#)(1))\times
  \operatorname{H}^1(\Gamma,(Y_2/Y_1^\#)(1))\to
  \operatorname{H}^2(\Gamma,\mathbb{Z}/n(2))\cong\mathbb{Z}/n$ 
  is a pairing of groups, i.e., a non-degenerate form.
\end{proof}

If two submodules $Y'_1\subset Y^\#_2$ and $Y'_2\subset Y^\#_1$ 
satisfy the equalities
$\operatorname{rank}Y'_i=\operatorname{rank}Y_i$ for $i=1,2$, 
then the inclusion $Y_i'\subset Y_i$ induces
the connecting homomorphism
$\partial_{Y_i/Y_i'}\colon T_i=(Y_i\otimes\overline F^\times)^\Gamma
\to\operatorname H^1(\Gamma,(Y_i/Y_i')(1))$.
Thus, we may compose the homomorphisms 
$\partial_{Y_1/Y_1'}$ and $\partial_{Y_2/Y_2'}$ with the cup product
to obtain a map
\[
  \operatorname{comm}_B\colon
  T_1\times T_2\stackrel{\partial\times\partial}\to
  \operatorname H^1(\Gamma,(Y_1/Y_1')(1))\times
  \operatorname H^1(\Gamma,(Y_2/Y_2')(1))\stackrel B\to
  \operatorname H^2(\Gamma,(\mathbb Z/n)(2)).
\]

\begin{lmm}
  \label{tCmm}
  For $i=1,2,3,4$, let $T_i$ be an algebraic torus 
  with cocharacter lattice $Y_i$, and
  $Y_i'\subset Y_i$ a $\Gamma$-submodule 
  of the same rank as $Y_i$.
  For $i=1,3$,
  we assume the inclusions $Y_i'\subset Y_{i+1}^\#$ and 
  $Y_{i+1}'\subset Y_{i}^\#$,
  so that a $\Gamma$-invariant bilinear form
  $B_{i,i+1}\colon Y_i\times Y_{i+1}\to\mathbb Z$ 
  induces a bilinear form
  $Y_i/Y_i'\times Y_{i+1}/Y'_{i+1}\to\mathbb Z/n$.
  For $i=1,2$, let $f_i\colon Y_i\to Y_{i+2}$ be 
  a $\Gamma$-equivariant homomorphism such that
  $f_i(Y_i')\subset Y_{i+2}'$ and that the diagram
  \[
    \begin{array}{*9c}
      Y_1/Y_1'\times Y_2/Y_2' &\stackrel{B_{1,2}}\to& \mathbb Z/n\\
      \mbox{}\downarrow\rlap{$\scriptstyle f_1\times f_2$} && 
      \rB{90}=\\
      Y_3/Y_3'\times Y_4/Y_4' &\stackrel{B_{3,4}}\to& \mathbb Z/n
    \end{array}
  \]
  commutes.
  Then these maps induce the following commutative diagram:
  \[
    \begin{array}{*9c}
      T_1\times T_2
      &\stackrel{\partial\times\partial}\to&
      \operatorname H^1(\Gamma,(Y_1/Y_1')(1))\times
      \operatorname H^1(\Gamma,(Y_2/Y_2')(1))
      &\stackrel{B_{1,2}}\to&
      \operatorname H^2(\Gamma,(\mathbb Z/n)(2))\\
      \mbox{}\downarrow\rlap{$\scriptstyle f_1\times f_2$} && 
      \downarrow &&
      \rB{90}=\\
      T_3\times T_4
      &\stackrel{\partial\times\partial}\to&
      \operatorname H^1(\Gamma,(Y_3/Y_3')(1))\times
      \operatorname H^1(\Gamma,(Y_4/Y_4')(1))
      &\stackrel{B_{3,4}}\to&
      \operatorname H^2(\Gamma,(\mathbb Z/n)(2)).
    \end{array}
  \]
\end{lmm}
\begin{proof}
  This follows from the functoriality of cup products.
\end{proof}

\begin{prp}
  Let $T_i$, $Y_i\supset Y_i'$ and $B_{i,i+1}$ be as above.
  For $i=1,2$, let $f_i\colon Y_i\to Y_{i+2}$ be 
  any $\Gamma$-equivariant homomorphism making the diagram
  \[
    \begin{array}{*9c}
      Y_1\times Y_2 &\stackrel{B_{1,2}}\to& \mathbb Z\\
      \mbox{}\downarrow\rlap{$\scriptstyle f_1\times f_2$} && 
      \rB{90}=\\
      Y_3\times Y_4 &\stackrel{B_{3,4}}\to& \mathbb Z
    \end{array}
  \]
  commute.
  Then the following diagram commutes:
  \[
    \begin{array}{*9c}
      T_1\times T_2
      &\stackrel{\partial\times\partial}\to&
      \operatorname H^1(\Gamma,(Y_1/Y_1')(1))\times
      \operatorname H^1(\Gamma,(Y_2/Y_2')(1))
      &\stackrel{B_{1,2}}\to&
      \operatorname H^2(\Gamma,(\mathbb Z/n)(2))\\
      \mbox{}\downarrow\rlap{$\scriptstyle f_1\times f_2$} && 
      &&
      \rB{90}=\\
      T_3\times T_4
      &\stackrel{\partial\times\partial}\to&
      \operatorname H^1(\Gamma,(Y_3/Y_3')(1))\times
      \operatorname H^1(\Gamma,(Y_4/Y_4')(1))
      &\stackrel{B_{3,4}}\to&
      \operatorname H^2(\Gamma,(\mathbb Z/n)(2)).
    \end{array}
  \]
  Especially, 
  the composite map 
  $\operatorname{comm}_{B_{1,2}}=B_{1,2}\circ(\partial\times\partial)$
  is independent of the choice of subgroups $Y_1'$ and $Y_2'$.
\end{prp}

\begin{proof}
  \begin{enumerate}[wide]
    \item\label{eBud}
      First, suppose that $Y_3'=nY_3$ and that $Y_4'=nY_4$.
      Since $f_i(nY_i)\subset nY_{i+2}$ for $i=1,2$,
      we have the following commutative diagram by \Cref{tCmm}:
      \[
        \begin{array}{*9c}
          T_1\times T_2
          &
          \to&
          \operatorname H^1(\Gamma,(Y_1/Y_1')(1))\times
          \operatorname H^1(\Gamma,(Y_2/Y_2')(1))
          &\stackrel{B_{1,2}}\to&
          \operatorname H^2(\Gamma,(\mathbb Z/n)(2))\\
          \rB{90}=&& 
          \downarrow &&
          \rB{90}=\\
          T_1\times T_2
          &
          \to&
          \operatorname H^1(\Gamma,(Y_1/Y_2^\#)(1))\times
          \operatorname H^1(\Gamma,(Y_2/Y_1^\#)(1))
          &\stackrel{B_{1,2}}\to&
          \operatorname H^2(\Gamma,(\mathbb Z/n)(2))\\
          \rB{90}=&& 
          \uparrow &&
          \rB{90}=\\
          T_1\times T_2
          &
          \to&
          \operatorname H^1(\Gamma,(Y_1/nY_1)(1))\times
          \operatorname H^1(\Gamma,(Y_2/nY_2)(1))
          &\stackrel{B_{1,2}}\to&
          \operatorname H^2(\Gamma,(\mathbb Z/n)(2))\\
          \mbox{}\downarrow\rlap{$\scriptstyle f_1\times f_2$} && 
          \downarrow &&
          \rB{90}=\\
          T_3\times T_4
          &
          \to&
          \operatorname H^1(\Gamma,(Y_3/nY_3)(1))\times
          \operatorname H^1(\Gamma,(Y_4/nY_4)(1))
          &\stackrel{B_{3,4}}\to&
          \operatorname H^2(\Gamma,(\mathbb Z/n)(2)).
        \end{array}
      \]

    \item
      We may reduce the general case to \ref{eBud} by the diagram
      \[
        \begin{array}{*9c}
          T_1\times T_2
          &\to&
          \operatorname H^1(\Gamma,(Y_1/Y_1')(1))\times
          \operatorname H^1(\Gamma,(Y_2/Y_2')(1))
          &\stackrel{B_{1,2}}\to&
          \operatorname H^2(\Gamma,(\mathbb Z/n)(2))\\
          \mbox{}\downarrow\rlap{$\scriptstyle f_1\times f_2$} && 
          &&
          \rB{90}=\\
          T_3\times T_4
          &\to&
          \operatorname H^1(\Gamma,(Y_3/nY_3)(1))\times
          \operatorname H^1(\Gamma,(Y_4/nY_4)(1))
          &\stackrel{B_{3,4}}\to&
          \operatorname H^2(\Gamma,(\mathbb Z/n)(2))\\
          \mbox{}\uparrow\rlap{$\scriptstyle\operatorname{id}\times
          \operatorname{id}$} && 
          &&
          \rB{90}=\\
          T_3\times T_4
          &\to&
          \operatorname H^1(\Gamma,(Y_3/Y_3')(1))\times
          \operatorname H^1(\Gamma,(Y_4/Y_4')(1))
          &\stackrel{B_{3,4}}\to&
          \operatorname H^2(\Gamma,(\mathbb Z/n)(2)).
        \end{array}
      \]
  \end{enumerate}
\end{proof}

\begin{crr}
  \label{tCmr}
  For $i=1,2$, let $Y_i\subset Y$ be a $\Gamma$-submodule 
  such that 
  the associated homomorphism $T_i\to T$ is injective. 
  Then the restriction 
  $B_{\widetilde T}|_{Y_1\times Y_2}\colon Y_1\times Y_2\to\mathbb Z$
  gives the equality
  \[
    \operatorname{comm}_{B_{\widetilde T}|_{Y_1\times Y_2}}=
    \operatorname{comm}_{\widetilde T}|_{T_1\times T_2}
  \]
  as bilinear forms 
  $T_1\times T_2\to\operatorname H^2(\Gamma,(\mathbb Z/n)(2))$.
\end{crr}

Note that
$\operatorname{comm}_{B_{\widetilde T}}=
\operatorname{comm}_{\widetilde T}$ for $B=B_{\widetilde T}$.

\begin{dfn}[annihilators]
  \label{tAnn}
  \begin{enumerate}
    \item
      Generally,
      let $A_1$ and $A_2$ be finite $\Gamma$-modules, and
      $B\colon A_1\times A_2\to(\mathbb Z/n)(2)$ 
      a $\Gamma$-equivariant bilinear form. 
      For a subgroup $J_1\subset\operatorname{H}^1(\Gamma,A_1)$,
      we define the subgroup 
      $J_1^\perp=J_1^{\perp_B}$ of $\operatorname{H}^1(\Gamma,A_2)$
      as the annihilator for the bilinear form
      \[
        B(\;\cup\;)\colon
        \operatorname{H}^1(\Gamma,A_1)\times
        \operatorname{H}^1(\Gamma,A_2)
        \to\operatorname{H}^2(\Gamma,\mathbb{Z}/n(2)).
      \]
      That is:
      \[
        J_1^\perp=\{h\in\operatorname{H}^1(\Gamma,A_2)\mid
        \forall j_1\in J_1, B(j_1\cup h)=0\}.
      \]
      Similarly,
      for a subgroup $J_2\subset\operatorname{H}^1(\Gamma,A_2)$,
      we define the subgroup 
      $J_2^\perp=J_2^{\perp_B}$ of $\operatorname{H}^1(\Gamma,A_1)$ by
      \[
        J_2^\perp=\{h\in\operatorname{H}^1(\Gamma,A_1)\mid
        \forall j_2\in J_2, B(h\cup j_2)=0\}.
      \]

    \item\label{eZdg}
      For subgroups $T_1,T_2\subset T$,
      we define the subgroup 
      $\operatorname Z^\dagger_{T_1}(T_2)
      \subset T_1$ 
      as the annihilator of $T_2$ for the commutator map
      $\operatorname{comm}_{\widetilde T}\colon T\times T\to\mu_n$.
      That is:
      \[
        \operatorname Z^\dagger_{T_1}(T_2)=
        \{t\in T_1\mid \operatorname{comm}_{\widetilde T}(t,T_2)=
        \{1\}\text{ in }\mu_n\}.
      \]
      Note that for $T_2=T$, 
      we have
      $\operatorname{Z}^\dagger_{T_1}(T)=\mathrm{Z}^\dagger\cap T_1$.
  \end{enumerate}
\end{dfn}

\begin{exm}
  \label{trAn}
  \begin{enumerate}
    \item\label{etAn}
      For $i=1,2$, suppose that 
      $T_i\subset T$ is a subtorus with cocharacter lattice 
      $Y_i\subset Y$.
      Then the restriction 
      $B_{\widetilde T}|_{Y_1\times Y_2}\colon 
      Y_1\times Y_2\to\mathbb Z$
      also describes the annihilator of a subgroup $S\subset T_2$ as
      \[
        \operatorname Z^\dagger_{T_1}(S)=\{t\in T_1\mid 
        \operatorname{comm}_{B_{\widetilde T}|Y_1\times Y_2}(t,S)=
        \{1\}\text{ in }\mu_n\}
      \]
      by \Cref{tCmr}.
      Take $\Gamma$-submodules $Y'_1\subset Y_2^\#$ and
      $Y'_2\subset Y_1^\#$ 
      such that $\operatorname{rank}Y_i'=\operatorname{rank}Y_i$
      for $i=1,2$.
      Then, we may write
      \begin{align*}
        \operatorname Z^\dagger_{T_1}(S)
        &=\{t\in T_1\mid \forall s\in S,
        B(\partial_{Y_1/Y_1'}(t)\cup\partial_{Y_2/Y_2'}(s))=0\}\\
        &=\partial_{Y_1/Y_1'}^{-1}(\partial_{Y_2/Y_2'}(S)^\perp),
      \end{align*}
      where the annihilator $\partial_{Y_2/Y_2'}(S)^\perp$ is 
      defined for the bilinear form
      $\operatorname H^1(\Gamma,(Y_1/Y_1')(1))\times
      \operatorname H^1(\Gamma,(Y_2/Y_2')(1))
      \stackrel{B}\to\operatorname H^2(\Gamma,(\mathbb Z/n)(2))$.

    \item
      \label{tDp}
      For the connecting homomorphism
      $\partial=\partial_{Y/Y^\#}\colon
      T=(Y\otimes\overline F^\times)^\Gamma\to
      \operatorname H^1(\Gamma,(Y/Y^\#)(1))$,
      we shall show in \Cref{Weiss} that
      the annihilator of the subgroup 
      $\partial((Y\otimes\overline{\mathcal{O}}^\times)^\Gamma)
      \subset\operatorname H^1(\Gamma,(Y/Y^\#)(1))$
      is itself under the same assumption as our main theorem:
      \[
        \partial((Y\otimes\overline{\mathcal{O}}^\times)^\Gamma)^\perp
        =\partial((Y\otimes\overline{\mathcal{O}}^\times)^\Gamma).
      \]
      By \ref{etAn}, we have
      \begin{align*}
        \operatorname Z^\dagger_T
        ((Y\otimes\overline{\mathcal O}^\times)^\Gamma)=
        \partial^{-1}\partial
        ((Y\otimes\overline{\mathcal O}^\times)^\Gamma)
        &=(Y\otimes\overline{\mathcal O}^\times)^\Gamma+
        \operatorname{Ker}\partial_{Y/Y^\#}\\
        &=(Y\otimes\overline{\mathcal O}^\times)^\Gamma+\iota(T^\#).
      \end{align*}
  \end{enumerate}
\end{exm}

The following is another example of the annihilator 
$\operatorname Z^\dagger_{T_1}(T_2)$
that can be described explicitly: 
\begin{prp}
  \label{tStz}
  For $i=1,2$, let $T_i\subset T$ be a subtorus
  with cocharacter lattice $Y_i\subset Y$.
  Recall that the 
  bilinear form
  $B_{\widetilde T}|_{Y_1\times Y_2}\colon Y_1\times Y_2\to\mathbb Z$
  defines the subgroups $Y_2^\#\subset Y_1$ and $Y_1^\#\subset Y_2$,
  which induce isogenies
  $T_2^\#\stackrel\iota\to T_1$ and $T_1^\#\stackrel\iota\to T_2$
  of tori (\Cref{dMT}).
  Assume that 
  the torus $T_1^\#$ splits over $F$, i.e. 
  $\Gamma$ acts trivially on the submodule $Y_1^\#\subset Y_2$.
  Then,
  \[
    \operatorname Z^\dagger_{T_1}(T_2)=
    \operatorname{Ker}\partial_{Y_1/Y_2^\#}=\iota(T^\#_2).
  \]
\end{prp}

\begin{proof}
  The last equality follows from the long exact sequence
  where the connecting homomorphism $\partial_{Y_1/Y_2^\#}$ appears.
  Thus, it suffices to see the first equality.
  Since
  $\operatorname Z^\dagger_{T_1}(T_2)
  =\partial_{Y_1/Y_2^\#}^{-1}(\partial_{Y_2/Y_1^\#}(T_2)^\perp)$
  by \Cref{trAn}~\ref{etAn}, it suffices to see that
  $\partial_{Y_2/Y_1^\#}(T_2)^\perp=0$.
  By assumption, the connecting homomorphism
  $\partial_{Y_2/Y_1^\#}\colon T_2=
  (Y_2\otimes\overline F^\times)^\Gamma\to
  \operatorname H^1(\Gamma,(Y_2/Y_1^\#)(1))$
  is surjective (\Cref{ePbi}~\ref{eDsj}).
  Since the bilinear form 
  \[
    \operatorname{H}^1(\Gamma,(Y_1/Y_2^\#)(1))\times
    \operatorname{H}^1(\Gamma,(Y_2/Y_1^\#)(1))
    \stackrel B\to\operatorname{H}^2(\Gamma,\mathbb{Z}/n(2))
  \]
  is non-degenerate by \Cref{eNdg}, we have
  \[
    \partial_{Y_2/Y_1^\#}(T_2)^\perp=
    \operatorname H^1(\Gamma,(Y_2/Y_1^\#)(1))^\perp=0.
  \]
\end{proof}

\begin{exm}
  \label{eGim}
  The maximal split subtorus $T^\Gamma\subset T$ is expressed as
  $T^\Gamma=Y^\Gamma\otimes F^\times$.
  Note that the exact sequence 
  $0\to Y^\Gamma\hookrightarrow Y\to Y/Y^\Gamma\to 0$ splits,
  since the quotient $Y/Y^\Gamma$ is a free abelian group.
  Hence the sequence
  \[
    0\to
    \begin{array}[t]{c}
      (Y^\Gamma\otimes\overline F^\times)^\Gamma\\
      \rB{90}=\\
      T^\Gamma
    \end{array}
    \to
    \begin{array}[t]{c}
      (Y\otimes\overline F^\times)^\Gamma\\
      \rB{90}=\\
      T
    \end{array}
    \to
    ((Y/Y^\Gamma)\otimes\overline F^\times)^\Gamma\to 0
  \]
  is split exact.
  Since the action of $\Gamma$ on $Y^\Gamma$ is trivial,
  the pair $(T,T^\Gamma)$ satisfies the assumption of \Cref{tStz}.
  Therefore
  \[\operatorname Z^\dagger_{T}(T^\Gamma)=\iota(T^{\Gamma\#}),\]
  where $T^{\Gamma\#}\stackrel\iota\to T$ is the isogeny
  corresponding to the inclusion $Y^{\Gamma\#}\hookrightarrow Y$
  of cocharacter lattices.
\end{exm}

\section{Main theorem}
\label{Weiss}
Our main theorem aims to determine the finite abelian group 
$\mathcal S_{\widetilde T}=\mathrm Z^\dagger/\iota(T^\#)$ introduced 
in \Cref{Das}.
Recall that its dual group 
$\operatorname{Hom}(\mathcal S_{\widetilde T},\mathbb C^\times)$ 
parametrizes a packet in the local Langlands correspondence
for a Brylinski-Deligne covering $\mu_n\to\widetilde T\to T$ of 
a torus over a non-archimedean local field $F$ of characteristic zero.
Here, $\mu_n\subset \mathbf f^\times$ denotes the cyclic subgroup
of order $n$.

We prove the main theorem 
by reducing it to a comparison between the image
$\mathrm{Z}^\dagger\subset T$
of the center $\operatorname Z(\widetilde T)\subset\widetilde T$
and the image of the isogeny $T^\#\stackrel\iota\to T$ of tori.
We analyze these two subgroups of $T$ via the exact sequence
$0\to(Y\otimes\overline{\mathcal O}^\times)^\Gamma\to T
\stackrel{\operatorname{ord}}\to Y^\Gamma\otimes\mathbb Q$
appeared in \Cref{ofq_Gal}.
In addition to the cocharacter lattices $Y^\#\subset Y$
of $T^\#$ and $T$, we employ the subgroup 
\[
  Y^{\Gamma\#}=(Y^\Gamma)^\#=
  \{y\in Y\mid B_{\widetilde T}(y,Y^\Gamma)\subset n\mathbb Z\}
\]
of $Y$, 
where $B_{\widetilde T}\colon Y\times Y\to\mathbb Z$  
is the bilinear form 
attached to the Brylinski-Deligne covering $\widetilde T\to T$.
Then the inclusion $Y^{\Gamma\#}\hookrightarrow Y$ again induces
an isogeny $T^{\Gamma\#}\stackrel\iota\to T$ of tori.

The essential part of the proof is
the self-orthogonality of the subgroup
$\partial_{Y/Y^\#}((Y\otimes\overline{\mathcal{O}}^\times)^\Gamma)
\subset\operatorname H^1(\Gamma,(Y/Y^\#)(1))$
mentioned in \Cref{trAn}~\ref{tDp},
which will be reduced to the orthogonality theorem 
for unramified cohomology.
Let $I=I_F\subset\Gamma$ be the inertia group of the extension
$\overline F/F$, and 
$\Gamma_{\mathbf f}=
\operatorname{Gal}(\overline{\mathbf f}/\mathbf f)=\Gamma/I$
the absolute Galois group of $\mathbf f$.
We fix a finite Galois extension $L/F$ of the base field where
the torus $T=(Y\otimes\overline F^\times)^\Gamma$ splits over $L$,
and write $e\in\mathbb Z_{>0}$ for its ramification index.

By the assumption
$\mu_n\subset\mathbf f^\times$,
we can employ the exact sequence
$0\to (Y/Y^\#)(1)\to Y^\#\otimes\overline{\mathbf f}^\times\to
Y\otimes\overline{\mathbf f}^\times\to 0$ given in \Cref{gSf}.

\begin{lmm}
  \label{tsI}
  Assume that the degree $n$ of the cover is prime to 
  the ramification index $e$ of the fixed splitting field $L/F$.
  \begin{enumerate}
    \item \label{gSfi}
      The inclusion $Y^\#\hookrightarrow Y$ induces 
      a short exact sequence
      \[
        0\to(Y/Y^\#)^I(1)\to(Y^\#\otimes\overline{\mathbf f}^\times)^I
        \to(Y\otimes\overline{\mathbf f}^\times)^I\to 0.
      \]
    \item\label{efI}
      The exact sequence in \ref{gSfi} induces an isomorphism
      \[
        \iota\colon\operatorname H^1(\Gamma_{\mathbf f},
        (Y^\#\otimes\overline{\mathbf f}^\times)^I)\to
        \operatorname H^1(\Gamma_{\mathbf f},
        (Y\otimes\overline{\mathbf f}^\times)^I).
      \]
    \item \label{tSd}
      The connecting homomorphism
      \[
        \partial\colon((Y\otimes\overline{\mathbf f}^\times)^I)
        ^{\Gamma_{\mathbf f}}\to
        \operatorname H^1(\Gamma_{\mathbf f},(Y/Y^\#)^I(1))
      \]
      induced from the exact sequence in \ref{gSfi} is surjective.
  \end{enumerate}
\end{lmm}

\begin{proof}
  \ref{gSfi}
  Let $I_L\subset\operatorname{Gal}(\overline F/L)$ be 
  the inertia group 
  for $L$.
  Note that $I_L$ acts trivially on 
  $Y\otimes\overline{\mathbf f}^\times$ and
  $Y^\#\otimes\overline{\mathbf f}^\times$, 
  and that $\#(I/I_L)=e$.
  Then the sequence in \Cref{gSf}~\ref{eSf} gives the exact sequence
  \[
    (Y^\#\otimes\overline{\mathbf f}^\times)^{I/I_L}\to
    (Y\otimes\overline{\mathbf f}^\times)^{I/I_L}\to
    \operatorname H^1(I/I_L,(Y/Y^\#)(1)).
  \]
  The last term of this sequence vanishes,
  since $(Y/Y^\#)(1)$ is a $\mathbb Z/n$-module 
  and the integers $n$ and $e$ are relatively prime by 
  assumption~\cite[Proposition 6.1.10]{Weib_Homlg}.

  \ref{efI}
  Since the cohomological dimension of the finite field $\mathbf f$
  is one, the two maps 
  \[
    \operatorname H^1(\Gamma_{\mathbf f},
    (Y\otimes\overline{\mathbf f}^\times)^I)\stackrel n\to
    \operatorname H^1(\Gamma_{\mathbf f},
    (Y^\#\otimes\overline{\mathbf f}^\times)^I)\stackrel\iota\to
    \operatorname H^1(\Gamma_{\mathbf f},
    (Y\otimes\overline{\mathbf f}^\times)^I)
  \]
  are surjective.
  Since $(Y\otimes\overline{\mathbf f}^\times)^I$ is
  an affine algebraic group defined over $\mathbf f$, we know that
  the group 
  $\operatorname H^1(\Gamma_{\mathbf f},
  (Y\otimes\overline{\mathbf f}^\times)^I)$
  has finite order~\cite[III.4.3]{Srr_Gal_coh}.
  Hence the composite 
  $\operatorname H^1(\Gamma_{\mathbf f},
  (Y\otimes\overline{\mathbf f}^\times)^I)\stackrel n\to
  \operatorname H^1(\Gamma_{\mathbf f},
  (Y\otimes\overline{\mathbf f}^\times)^I)$
  is an isomorphism, and so is the map $\iota$.

  \ref{tSd}
  In the exact sequence
  \[
    ((Y\otimes\overline{\mathbf f}^\times)^I)^{\Gamma_{\mathbf f}}
    \stackrel\partial\to
    \operatorname H^1(\Gamma_{\mathbf f},(Y/Y^\#)^I(1))\to
    \operatorname H^1(\Gamma_{\mathbf f},
    (Y^\#\otimes\overline{\mathbf f}^\times)^I)\stackrel\iota\to
    \operatorname H^1(\Gamma_{\mathbf f},
    (Y\otimes\overline{\mathbf f}^\times)^I),
  \]
  the middle map is zero, since the last map $\iota$ is 
  an isomorphism.
\end{proof}

For a $\Gamma$-module $A$, let
$\inf\colon\operatorname H^1(\Gamma_{\mathbf f},A^I)\to
\operatorname H^1(\Gamma,A)$
be the inflation map on the first cohomology group.
This $\inf$ is injective, and fits into the exact sequence
in low dimensions given by the Lyndon-Hochschild-Serre 
spectral sequence.
Thus, we often regard
$\operatorname H^1(\Gamma_{\mathbf f},A^I)\subset
\operatorname H^1(\Gamma,A)$.

\begin{prp}
  \label{eUrc}
  Let $A\times A'\stackrel B\to\mu_n$ be
  a non-degenerate $\Gamma$-equivariant pairing between two
  finite $\Gamma$-modules.
  Then, in the non-degenerate pairing 
  $\beta\colon\operatorname H^1(\Gamma,A)\times
  \operatorname H^1(\Gamma,A')\stackrel{B(\,\cup\,)}\to
  \operatorname H^2(\Gamma,\mu_n)=\mathbb Z/n$ due to the local 
  Tate duality, we have
  \[
    \inf(\operatorname{H}^1(\Gamma_{\mathbf{f}},A^I))^\perp=
    \inf(\operatorname{H}^1(\Gamma_{\mathbf{f}},A'^I)).
  \]
\end{prp}

\begin{proof}
  The restriction 
  $\operatorname H^1(\Gamma_{\mathbf f},A^I)\times
  \operatorname H^1(\Gamma_{\mathbf f},A'^I)\to
  \operatorname H^2(\Gamma,\mu_n)$
  of $\beta$ 
  factors through $\operatorname H^2(\Gamma_{\mathbf f},\mu_n)=0$.
  That is,
  \[
    \beta\left(\operatorname H^1(\Gamma_{\mathbf f},A^I),\,
    \operatorname H^1(\Gamma_{\mathbf f},A'^I)\right)=0.
  \]
  It remains to show that
  $\#\operatorname H^1(\Gamma,A)=
  \#\operatorname H^1(\Gamma_{\mathbf f},A^I)\cdot
  \#\operatorname H^1(\Gamma_{\mathbf f},A'^I)$.
  Indeed, we see the equalities
  \begin{align*}
    \#\operatorname H^1(\Gamma,A)
    &=\#\operatorname H^0(\Gamma,A)\cdot\#\operatorname H^2(\Gamma,A)
    \\
    &=\#\operatorname H^0(\Gamma,A)\cdot\#\operatorname H^0(\Gamma,A')
    \\
    &=\#\operatorname H^0(\Gamma_{\mathbf f},A^I)\cdot
    \#\operatorname H^0(\Gamma_{\mathbf f},A'^I) \\
    &=\#\operatorname H^1(\Gamma_{\mathbf f},A^I)\cdot
    \#\operatorname H^1(\Gamma_{\mathbf f},A'^I)
  \end{align*}
  as follows:

  For the first equality, note that $A$ is a finite 
  $\mathbb Z/n$-module by definition.
  Since $n$ is prime to the residual characteristic $p$
  by assumption, the Euler-Poincar\'e characteristic is
  $\frac{\#\operatorname H^0(\Gamma,A)\cdot
  \#\operatorname H^2(\Gamma,A)}{\#\operatorname H^1(\Gamma,A)}=1$.
  The second equality follows from the local Tate duality
  between $\operatorname H^2$ and $\operatorname H^0$.
  The third one holds by definition.
  For an explicit proof of the last equality, we may consult
  Harari's textbook~\cite[Lemma 10.14]{Hara_Gal}.
\end{proof}

Note that this lemma slightly generalizes 
the well-known orthogonality theorem for unramified 
cohomology~\cite[II.5.5]{Srr_Gal_coh}.

\begin{crr}
  \label{tOd}
  Assume that the degree $n$ of the cover is prime to 
  the ramification index $e$ of the fixed splitting field $L/F$.
  \begin{enumerate}
    \item\label{eOdd}
      For the connecting homomorphism
      $\partial=\partial_{Y/Y^\#}\colon
      T=(Y\otimes\overline F^\times)^\Gamma\to
      \operatorname H^1(\Gamma,(Y/Y^\#)(1))$,
      the annihilator of the subgroup 
      $\partial((Y\otimes\overline{\mathcal{O}}^\times)^\Gamma)
      \subset\operatorname H^1(\Gamma,(Y/Y^\#)(1))$
      is itself:
      \[
        \partial((Y\otimes\overline{\mathcal{O}}^\times)^\Gamma)
        ^\perp
        =\partial((Y\otimes\overline{\mathcal{O}}^\times)^\Gamma).
      \]
    \item\label{eZo}
      As subgroups of $T=(Y\otimes\overline F^\times)^\Gamma$,
      we have 
      \[
        \mathrm Z^\dagger\subset
        \operatorname Z^\dagger_T
        ((Y\otimes\overline{\mathcal O}^\times)^\Gamma)=
        (Y\otimes\overline{\mathcal O}^\times)^\Gamma+\iota(T^\#).
      \]
      Especially, the map
      \[
        \left.\mathrm Z^\dagger\cap
        (Y\otimes\overline{\mathcal O}^\times)^\Gamma\middle/
        \iota(T^\#)\cap
        (Y\otimes\overline{\mathcal O}^\times)^\Gamma\right.\to
        \mathrm Z^\dagger/\iota(T^\#)=\mathcal S_{\widetilde T}
      \]
      is an isomorphism.
  \end{enumerate}
\end{crr}

\begin{proof}
  \ref{eOdd}
  By assumption,
  \Cref{tsI}~\ref{gSfi} ensures that the horizontal sequences in
  the commutative diagram
  \[
    \begin{array}{*9c}
      0&\to& (Y/Y^\#)(1)  &\to& Y^\#\otimes\overline{\mathcal O}^\times  &\to& Y\otimes\overline{\mathcal O}^\times &\to& 0\\
      &   & \rB{90}=     &   & \downarrow                               &   & \rB{-90}{$\twoheadrightarrow$}\rlap{\ $\scriptstyle\text{split}$}                           &   &\\
      0&\to& (Y/Y^\#)(1)  &\to& Y^\#\otimes\overline{\mathbf f}^\times   &\to& Y\otimes\overline{\mathbf f}^\times  &\to& 0\\
      &   & \uparrow     &   & \uparrow                                 &   & \uparrow                             &   &\\
      0&\to& (Y/Y^\#)^I(1)&\to& (Y^{\#}\otimes\overline{\mathbf f}^\times)^{I}&\to& (Y\otimes\overline{\mathbf f}^\times)^I&\to& 0
    \end{array}
  \]
  are exact.
  The connecting homomorphisms $\partial$ on
  the Galois cohomology groups form 
  the commutative diagram
  \[
    \begin{array}{*9c}
      (Y\otimes\overline{\mathcal O}^\times)^\Gamma               &\stackrel\partial\to&\operatorname H^1(\Gamma,(Y/Y^\#)(1))\\
      \rB{-90}{$\twoheadrightarrow$}                                                  &                    &\rB{90}=\\
      (Y\otimes\overline{\mathbf f}^\times)^\Gamma                &\stackrel\partial\to&\operatorname H^1(\Gamma,(Y/Y^\#)(1))\\
      \rB{90}=                                                    &                    &\rB{90}{$\hookrightarrow$}\rlap{$\scriptstyle\inf$}\\
      ((Y\otimes\overline{\mathbf f}^\times)^I)^{\Gamma_{\mathbf f}}&\stackrel\partial\twoheadrightarrow&\operatorname H^1(\Gamma_{\mathbf f},(Y/Y^\#)^I(1))\rlap.
    \end{array}
  \]
  In the left column, the downward map
  $(Y\otimes\overline{\mathcal O}^\times)^\Gamma\to
  (Y\otimes\overline{\mathbf f}^\times)^\Gamma$
  is surjective, since 
  $Y\otimes\overline{\mathcal O}^\times\to
  Y\otimes\overline{\mathbf f}^\times$
  is a split surjection.
  The bottom $\partial$ is also surjective by \Cref{tsI}~\ref{tSd}.
  Recall that $\inf$ in the right column is injective.
  Then
  $\partial((Y\otimes\overline{\mathcal{O}}^\times)^\Gamma)=
  \partial((Y\otimes\overline{\mathbf{f}}^\times)^\Gamma)=
  \inf\left(\operatorname{H}^1(\Gamma_{\mathbf{f}},
  (Y/Y^\#)^I(1))\right)$.
  By \Cref{eUrc}, 
  \[
    \partial((Y\otimes\overline{\mathcal{O}}^\times)^\Gamma)=
    \inf\left(\operatorname{H}^1(\Gamma_{\mathbf{f}},
    (Y/Y^\#)^I(1))\right)^\perp=
    \partial((Y\otimes\overline{\mathcal{O}}^\times)^\Gamma)^\perp.
  \]

  \ref{eZo}
  The inclusion holds by definition.
  By \Cref{trAn}~\ref{etAn}, we have
  \begin{align*}
    \operatorname Z^\dagger_T
    ((Y\otimes\overline{\mathcal O}^\times)^\Gamma)=
    \partial^{-1}\partial
    ((Y\otimes\overline{\mathcal O}^\times)^\Gamma)
    &=(Y\otimes\overline{\mathcal O}^\times)^\Gamma+
    \operatorname{Ker}\partial_{Y/Y^\#}\\
    &=(Y\otimes\overline{\mathcal O}^\times)^\Gamma+\iota(T^\#).
  \end{align*}
\end{proof}

\Cref{tOd} allows us to describe the intersection
$\mathrm Z^\dagger\cap(Y\otimes\overline{\mathcal O}^\times)^\Gamma$.

\begin{prp}
  \label{tZtg}
  As in \Cref{tOd},
  assume that the degree $n$ of the cover is prime to 
  the ramification index $e$ of the fixed splitting field $L/F$.
  Then the following equalities hold.
  \begin{gather}
    \setcounter{equation}{0}
    \label{qSGO}
    T=\iota(T^\#)+T^\Gamma+
    (Y\otimes\overline{\mathcal O}^\times)^\Gamma\\
    \label{eZtg}
    \mathrm Z^\dagger\cap(Y\otimes\overline{\mathcal O}^\times)^\Gamma
    =\iota((Y^{\Gamma\#}\otimes\overline{\mathcal O}^\times)^\Gamma).
  \end{gather}
\end{prp}

\begin{proof}
  \eqref{qSGO}
  As usual, we normalize the valuation map
  $\operatorname{ord}\colon\overline F^\times\to\mathbb Q$
  so that $\operatorname{ord}(\overline F^\times)=\mathbb Z$.
  Then, for the induced map
  $\operatorname{ord}\colon T=(Y\otimes\overline F^\times)^\Gamma\to
  (Y\otimes\mathbb Q)^\Gamma=Y^\Gamma\otimes\mathbb Q$, we have
  $\operatorname{ord}(T^\Gamma)=Y^\Gamma$ and
  $\operatorname{ord}(T)\subset\frac1eY^\Gamma$.
  Since $nY\subset Y^\#$ by definition, we see the inclusion 
  $n(Y\otimes\overline F)^\Gamma\subset
  \iota((Y^\#\otimes\overline F)^\Gamma)$
  in $T=(Y\otimes\overline F)^\Gamma$, and hence
  $n\cdot\operatorname{ord}(T)\subset\operatorname{ord}(\iota(T^\#))$.
  Since the integers $n$ and $e$ are relatively prime by assumption,
  \[
    \operatorname{ord}(T)=
    n\operatorname{ord}(T)+e\operatorname{ord}(T)\subset
    \operatorname{ord}(\iota(T^\#))+Y^\Gamma=
    \operatorname{ord}(\iota(T^\#)+T^\Gamma).
  \]
  Since the map 
  $\operatorname{ord}\colon T\to Y^\Gamma\otimes\mathbb Q$
  has kernel $(Y\otimes\overline{\mathcal O}^\times)^\Gamma$,
  the claimed equality holds.

  \eqref{eZtg}
  The identity in \eqref{qSGO} and \Cref{eGim} imply
  \begin{align*}
    \mathrm Z^\dagger
    &=\operatorname Z^\dagger_T(T^\Gamma)\cap
    \operatorname Z^\dagger_T\left(\iota(T^\#)+
    (Y\otimes\overline{\mathcal O}^\times)^\Gamma\right)\\
    &=\iota(T^{\Gamma\#})\cap
    \operatorname Z^\dagger_T\left(\iota(T^\#)+
    (Y\otimes\overline{\mathcal O}^\times)^\Gamma\right).
  \end{align*}
  Since 
  $\operatorname Z^\dagger_T\left(\iota(T^\#)+
  (Y\otimes\overline{\mathcal O}^\times)^\Gamma\right)\supset
  (Y\otimes\overline{\mathcal O}^\times)^\Gamma$
  by
  \Cref{tOd}~\ref{eZo}, we see
  $
  \mathrm Z^\dagger
  \cap(Y\otimes\overline{\mathcal O}^\times)^\Gamma
  =\iota(T^{\Gamma\#})
  \cap(Y\otimes\overline{\mathcal O}^\times)^\Gamma
  =\iota((Y^{\Gamma\#}\otimes\overline{\mathcal O}^\times)^\Gamma)
  $.
  Note that the last equality is the special case of 
  \Cref{ePbi}~\ref{eePbi} for $Y'=Y^{\Gamma\#}$.
\end{proof}

Now we can prove the main theorem:
\begin{thm}
  \label{tMt}
  Let $F$ be a non-archimedean local field of characteristic zero,
  $\mathbf f$ the residue field of the valuation ring of $F$, and
  $\mu_n\subset \mathbf f^\times$ the cyclic subgroup of order $n$.
  Suppose that an algebraic torus $T=\mathbb{T}(F)$ defined over $F$ 
  splits over a finite Galois extension $L/F$ 
  whose ramification index $e$ is relatively prime to $n$.
  Let $\mu_n\to\widetilde{T}\to T$ be a Brylinski-Deligne 
  covering group.

  Then the group $\mathcal S=\mathcal S_{\widetilde T}$ 
  parametrizing a packet of representations of $\widetilde{T}$ is 
  \begin{align*}
    \mathcal S_{\widetilde T}
    &=\left.\iota((Y^{\Gamma\#}\otimes_{\mathbb Z}
    \overline{\mathcal{O}}^\times)^\Gamma)\middle/
    \iota((Y^{\#}\otimes_{\mathbb Z}
    \overline{\mathcal{O}}^\times)^\Gamma)\right..\\
    &=\left.\iota((Y^{\Gamma\#}\otimes_{\mathbb Z}
    \overline{\mathbf f}^\times)^\Gamma)\middle/
    \iota((Y^{\#}\otimes_{\mathbb Z}
    \overline{\mathbf f}^\times)^\Gamma)\right..
  \end{align*}
\end{thm}

\begin{proof}
  The first identity follows from
  \Cref{tOd}~\ref{eZo}, \Cref{tZtg}~\eqref{eZtg}, and
  \Cref{ePbi}~\ref{eePbi} for $Y'=Y^\#$.

  Let $\overline{\mathfrak p}\subset\overline{\mathcal{O}}$ be
  the maximal ideal of the valuation ring in 
  the algebraic closure $\overline F$. 
  Since $\mu_n\subset\mathbf f^\times$ by assumption,
  the inclusions
  $nY\subset Y^\#\subset Y^{\Gamma\#}\subset Y$ induces
  isomorphisms
  \[
    Y^\#\otimes(1+\overline{\mathfrak p})\stackrel\sim\to
    Y^{\Gamma\#}\otimes(1+\overline{\mathfrak p})\stackrel\sim\to
    Y\otimes(1+\overline{\mathfrak p}).
  \]
  Then the split exact sequence
  $1\to1+\overline{\mathfrak p}\to
  \overline{\mathcal{O}}^\times
  \to\overline{\mathbf{f}}^\times\to1$
  gives a commutative diagram 
  \[
    \begin{matrix}
      0                                                &                     &0                                                        &                     &0                                            \\
      \downarrow                                       &                     &\downarrow                                               &                     &\downarrow                                   \\
      (Y^\#\otimes(1+\overline{\mathfrak p}))^\Gamma   &\stackrel{\sim}{\to} &(Y^{\Gamma\#}\otimes(1+\overline{\mathfrak p}))^\Gamma   &\stackrel{\sim}{\to} &(Y\otimes(1+\overline{\mathfrak p}))^\Gamma  \\
      \downarrow                                       &                     &\downarrow                                               &                     &\downarrow                                   \\
      (Y^\#\otimes\overline{\mathcal O}^\times)^\Gamma &\to                  &(Y^{\Gamma\#}\otimes\overline{\mathcal O}^\times)^\Gamma &\to                  &(Y\otimes\overline{\mathcal O}^\times)^\Gamma\\
      \downarrow                                       &                     &\downarrow                                               &                     &\downarrow                                   \\
      (Y^\#\otimes\overline{\mathbf f}^\times)^\Gamma  &\to                  &(Y^{\Gamma\#}\otimes\overline{\mathbf f}^\times)^\Gamma  &\to                  &(Y\otimes\overline{\mathbf f}^\times)^\Gamma \\
      \downarrow                                       &                     &\downarrow                                               &                     &\downarrow                                   \\
      0                                                &                     &0                                                        &                     &0\rlap{,}                                    \\
    \end{matrix}
  \]
  where the top horizontal maps are isomorphisms, and
  the vertical sequences are split exact.
  Thus the map
  \[
    \mathcal S_{\widetilde T}=
    \iota((Y^{\Gamma\#}\otimes\overline{\mathcal O}^\times)^\Gamma)
    /\iota((Y^\#\otimes\overline{\mathcal O}^\times)^\Gamma)\to
    \iota((Y^{\Gamma\#}\otimes\overline{\mathbf f}^\times)^\Gamma)
    /\iota((Y^\#\otimes\overline{\mathbf f}^\times)^\Gamma)
  \]
  is an isomorphism.
\end{proof}

The trivial example of the main theorem is
a Brylinski-Deligne cover $\widetilde T\to T$ of a split torus.
Then $\mathcal S_{\widetilde T}=1$, since $Y^\Gamma=Y$.
Another elementary example is a cover $\widetilde T\to T$ of 
an unramified torus, i.e., a torus
splitting over an unramified extension $L/F$ of the base field.
Then \Cref{tMt} regains Weissman's result~\cite{Weiss_loc_tori}
on covers of an unramified torus. 

Thus, our main theorem is a generalization of
Weissman's work on covering groups.
Moreover, this theorem partially realizes his 
hope~\cite[Remark 5.18]{Weiss_loc_tori} to parametrize a packet
in the local Langlands correspondence for
a cover of a ramified torus.

\end{document}